\documentclass[12pt, twoside, a4paper, fleqn]{article}

\usepackage{latexsym}%		Extra math symbols
\usepackage{amscd}%			AMS package for doing commutative diagrams
\usepackage{theorem}%		Improved theoremlike environments
\usepackage{pifont}%			Dingbat symbols (used for end of proof symbol)
\usepackage{mathbbol}%		Blackboard bold symbols
\usepackage{amsfonts}%		More fonts
\usepackage{xspace}%			Removes need for adding a \  after macros
\usepackage{amssymb}%		More math symbols
\usepackage{fancyhdr}%		Fancyheadings package
\usepackage{mathrsfs}%		For the 3rd letters
\usepackage{amsmath}%
\usepackage{graphics}%
\usepackage{graphicx}%
\usepackage{euscript}%
\usepackage{psfrag}%
\usepackage{empheq}%
\usepackage{upgreek}%		Upright Greek letters
\usepackage{eufrak}%			More fancy letters

\DeclareMathAlphabet{\mathpzc}{OT1}{pzc}{m}{it}%			More Math alphabets

{\theorembodyfont{\slshape} \newtheorem{theorem}{Theorem}}
{\theorembodyfont{\slshape} \newtheorem{lemma}[theorem]{Lemma}}
{\theorembodyfont{\slshape} \newtheorem{proposition}[theorem]{Proposition}}
{\theorembodyfont{\slshape} \newtheorem{corollary}[theorem]{Corollary}}

\newenvironment{proof}{\noindent\textbf{Proof:\ }}{$\hfill\Box$}

\numberwithin{equation}{section}

\title{\textsc{The Recurrence Rate and Hausdorff Dimension \\ of a Neighbourhood of some Typical Point \\ in the Julia Set of a Rational Map}}                %Full title

\author{Shrihari Sridharan \bigskip \\ Chennai Mathematical Institute (CMI), India. \bigskip \\ {\tt shrihari@cmi.ac.in}}     
\date{October 17, 2013}

\begin{document}

\maketitle

\bigskip
\begin{abstract}
\bigskip

\noindent 
In this article, we consider hyperbolic rational maps restricted on thier Julia sets and study about the recurrence rate of typical orbits in arbitrarily small neighbourhoods around them and their relationship to the Hausdorff dimension of such small neighbourhoods. 

\end{abstract}
\bigskip
\bigskip

\begin{tabular}{l l}
\textbf{Keywords} & Complex dynamics of rational maps, \\ & Recurrence rates, Hausdorff dimension. \bigskip \\
\textbf{AMS Subject Classifications} & 37F10, 37F15, 37F35.
\end{tabular}

\thispagestyle{empty}

\newpage

\section{Introduction}

Let $T : X \longrightarrow X$ be a transformation that preserves an ergodic probability measure $\mu$. Any neighborhood, say of radius $r > 0$, however small, around a generic point $x \in X$ in $X$, denoted by $N_{r} (x)$ must return to itself infinitely often, according to the Poincar\'{e} recurrence theorem. In this paper, we study the relationship between the Poincar\'{e} recurrence rate of typical points $y \in N_{r} (x) \subset X$ and the Hausdorff dimension of $N_{r} (x)$. We shall also focus on the case when $r \searrow 0$. In recent papers, Boshernitzan studied the relationship between recurrence rates and $\sigma$-finite outer measures for measure preserving dynamical systems focussing on billiards flows, symbolic systems and interval maps, \cite{mb:93}, Barreira and Saussol studied the quantitative behaviour of the recurrence rates by imposing a condition on the measure, \cite{bs:01} and Saussol studied the recurrence rates in rapidly mixing dynamical systems, \cite{bs:06}. In yet another recent paper, the author studied the relationship between incidence rate at which the forward orbit of a generic point $y \notin N_{r} (x)$ would reach $N_{r} (x)$ and the Hausdorff dimension of $N_{r} (x)$, \cite{ss:12}. This paper complements the study in \cite{ss:12}.  
\medskip 

\noindent 
In this paper, we prove analogous results in the setting of complex dynamics; hyperbolic rational maps restricted on their Julia sets and a non-atomic probability measure preserved by the rational map; for example, the Sinai - Ruelle - Bowen (SRB) measure.  
\medskip

\noindent 
In sections 2 and 3, we write the fundamental and involved definitions that constitute the skeleton on which this paper rests its studies on. We shall also encounter some simple properties, as much necessary of the terms defined therein. In section 4, we state the main results. In further sections, we prove the theorems stated in section 4. 
\bigskip

\section{Fundamental Definitions} 

Let $\overline{\mathbb{C}}$ denote the \emph{Riemann sphere} and let $T$ be a \emph{rational map} defined on the Riemann sphere. By \emph{degree} of the rational map, we mean the number of inverse images for a typical point $z \in \overline{\mathbb{C}}$ counted with multiplicity. In other words, the maximum among the degrees of the two relatively prime polynomials whose quotient yields the rational map is defined to be its degree denoted by $d$. For our purpose of study in this paper, we shall only consider those rational maps whose degree is at least $2$. One of the several possible definitions of the \emph{Julia set} $\mathcal{J} \subset \overline{\mathbb{C}}$ of $T$ states that it is the closure of the set of all repelling periodic points, i.e.,
\begin{equation} 
\mathcal{J}\ :=\ \overline{\left\{ z \in \overline{\mathbb{C}}\ :\ T^{p} z = z\ \ \text{for some}\ \ p \in \mathbb{Z}_{+}\ \ \text{and}\ \ |(T^{p})' (z)| > 1 \right\}}. 
\end{equation}
\medskip 

\noindent
Elementary observations reveal that the rational map remains completely invariant on its Julia set, i.e., $T^{-1} (\mathcal{J}) = \mathcal{J}$. For example, consider the polynomial map $T(z) := z^{d}$ defined on $\overline{\mathbb{C}}$. The Julia set of this polynomial map is the unit circle in the complex plane; 
\[ \mathcal{J}\ \ =\ \ \mathbb{S}^{1}\ \ :=\ \ \left\{ z \in \mathbb{C}\ :\ |z| = 1 \right\}. \] 
For more properties of Julia sets of rational maps, please refer \cite{ml:86}. We focus on \emph{hyperbolic rational maps} restricted on their Julia sets in this paper, i.e., there exists $C > 0$ and $\lambda > 1$ such that for all $z \in \mathcal{J}$ and $n \ge 1$, we have $|(T^{n})' (z)| \ge C \lambda^{n}$. Since the Julia set is compact and the transformation $T$ is continuous, the set of non-atomic $T$-invariant Borel probability measures defined on $\mathcal{J}$, denoted by $\mathcal{M}_{T} (\mathcal{J})$ is non-empty. The Lyubich's measure that equidistributes the pre-images of a typical point in $\mathcal{J}$ and the periodic points of $T$ in $\mathcal{J}$ is one such example. Observe that the Lyubich's measure reduces to the Haar measure on $\mathbb{S}^{1}$.
\medskip

\noindent 
We now define the \emph{pressure of a real-valued continuous function} in accordance with thermodynamic formalism. Consider a continuous function, $f : \mathcal{J} \longrightarrow \mathbb{R}$. Its pressure is defined by 
\begin{equation} 
\textrm{Pr} (f)\ \ :=\ \ \sup \left\{ h_{\mu} (T) + \int f d \mu\ \ :\ \ \mu \in \mathcal{M}_{T} (\mathcal{J}) \right\}. 
\end{equation} 
Here, $h_{\mu} (T)$ denotes the \emph{entropy} of the transformation with respect to the measure $\mu$, see \cite{pw:82} for more details. If $f$ is a real-valued, H\"{o}lder continuous function defined on the Julia set $\mathcal{J}$ of some hyperbolic rational map $T$, then, by a result due to Denker and Urbanski in \cite{du:91}, there exists a unique equilibrium measure, called the \emph{Sinai - Ruelle - Bowen measure} (the Gibbs' state) denoted by $\mu_{f} \in \mathcal{M}_{T} (\mathcal{J})$ realising the supremum in the definition of pressure. We further remark that by adding a coboundary to the H\"{o}lder continuous function $f$, one can normalise pressure, as done by Haydn in \cite{nh:99} so that $\textrm{Pr} (f) = 0$. 
\medskip

\noindent 
It is merely an observation that such systems $T : \mathcal{J} \longrightarrow \mathcal{J}$ alongwith $\mu \in \mathcal{M}_{T} (\mathcal{J})$ are ergodic. One of the several possible definitions of \emph{ergodicity} states, given any real-valued function $f \in L^{1} (\mu)$, 
\begin{equation} 
\lim_{n \to \infty} \frac{1}{n} \sum_{j = 0}^{n - 1} f \circ T^{j} (z)\ \ -\ \ \int f d \mu\ \ \ =\ \ \ 0,\ \ \ \ \mu \text{-a.e.} 
\end{equation}
Now let $f , g \in L^{2} (\mu)$. Consider the quantity called \emph{covariance} defined by 
\begin{equation} 
Cov_{T} (f, g)\ \ :=\ \ \lim_{n \to \infty} \int f \circ T^{n} g d \mu\ \ -\ \ \int f d \mu \int g d \mu. 
\end{equation} 
$T$ is said to be \emph{mixing} if $Cov_{T} (f, g) = 0$ for every $f, g \in L^{2} (\mu)$. 
\medskip 

\noindent 
In general, $Cov_{T} (f, g)$ vanishes at an arbitrarily slow pace as $n \to \infty$. It requires a little more structure on the underlying space for us to say anything tangible about the vanishing rate of $Cov_{T} (f, g)$. Fortunately, we are only dealing with Julia sets of rational maps in this paper and they have the requisite structure being compact metric spaces. An upper bound for $Cov_{T} (f, g)$ is then provided by $\| f \| \| g \| \theta_{n}$, where 
\[ \| f \|\ \ :=\ \ \sup_{z_{1}, z_{2} \in \mathcal{J}}\ \frac{| f (z_{1} - f (z_{2}) |}{| z_{1} - z_{2} |}, \]
and $\theta_{n} \to 0$, as $n \to \infty$. 
\medskip 

\noindent 
Dependent on the rate at which $\theta_{n}$ vanishes as $n \to \infty$, we shall call $Cov_{T} (f, g)$ to vanish as $n \to \infty$. For example, if $(\theta_{n})$ is a sequence of real numbers such that $n^{-p} \theta_{n} \to 0$ as $n \to \infty$ for some fixed $p > 0$, we shall say that $Cov_{T} (f, g)$ has a \emph{polynomial decay}. In this paper we shall focus on systems whose covariance have \emph{super-polynomial decay}, i.e., $n^{-p} \theta_{n} \to 0$ as $n \to \infty$ for every $p > 0$. 
\medskip 

\noindent 
Suppose we partition $\mathcal{J}$ finitely into $\{ N_{i} \}_{i = 1}^{M}$ such that $\bigcup \overline{N_{i}} = \mathcal{J}$ and $\mu ( N_{i} \cap N_{j} ) = 0$ whenever $i \ne j$. We call the decay as a \emph{local decay} (say on $N_{i}$) if $Cov_{T} (f, g) \le \| f \| \| g \| \theta_{i, n}$ where $supp (f) \subset N_{i},\ supp (g) \subset N_{i}$ and $\theta_{i, n}$ decays at some rate. Observe that if the covariance has a super-polynomial rate of decay, then the system has a super-polynomial rate of local decay in every component $N_{i}$ of the partition of $\mathcal{J}$. 
\bigskip 

\section{Recurrence Rate and Hausdorff Dimension}

Fix $z \in \mathcal{J}$ and for any $r > 0$, consider the neighbourhood $N_{r} (z)$ in $\mathcal{J}$ centered at $z$ of radius $r$. Owing to the non-atomicity of the measure $\mu \in \mathcal{M}_{T} (\mathcal{J})$, it is clear that $\mu (N_{r} (z)) > 0$ whenever $r > 0$. By the Poincar\'{e} recurrence theorem, we then have that the orbit of $\mu$-a.e. $w \in N_{r} (z)$ returns to $N_{r} (z)$ infinitely often. We define the Poincar\'{e} recurrence time of the centre $z \in N_{r} (z)$ as follows. 
\begin{equation} 
\tau_{r} (z)\ \ :=\ \ \inf \left\{n \in \mathbb{Z}_{+}\ :\ T^{n} (z) \in N_{r} (z) \right\}. 
\end{equation} 
More generally, allowing a slight abuse of notation, the Poincar\'{e} recurrence time of any typical point $w \in N_{r} (z)$ in $N_{r} (z)$ is denoted and defined by 
\begin{equation} 
\tau_{r} (w, z)\ \ :=\ \ \inf \left\{n \in \mathbb{Z}_{+}\ :\ T^{n} (w) \in N_{r} (z) \right\}. 
\end{equation} 
Thus, $\tau_{r} (z) \equiv \tau_{r} (z, z)$. In fact, this definition is valid even for $w \notin N_{r} (z)$. Hence, for any point $w \in \mathcal{J}$, one can define the incidence time of $w$ in $N_{r} (z)$ as 
\begin{equation} 
\tau_{r} (w, z)\ \ :=\ \ \inf \left\{n \in \mathbb{Z}_{+}\ :\ T^{n} (w) \in N_{r} (z) \right\}. 
\end{equation}  
There is a possibility that $T^{n} (w) \notin N_{r} (z),\ \forall n$. For example, consider a periodic point $w \in \mathcal{J}$ whose cycle remains away from $N_{r} (z)$. In such situations, we define 
\begin{equation} 
\tau_{r} (w, z)\ \ :=\ \ \infty.
\end{equation} 
Observe that in all these definitions, $\tau$\ and $r$\ are inversely proportional to each other. In particular, 
\begin{equation}
\tau_{kr} (z)\ \ \le\ \ \tau_{r} (z)\ \ \ \ \forall k \ge 1. 
\end{equation}
Similarly, 
\begin{equation}
\tau_{kr} (w, z)\ \ \le\ \ \tau_{r} (w, z)\ \ \ \ \forall w \in \mathcal{J}\ \ \ \text{and}\ \ \ \forall k \ge 1. 
\end{equation} 
Further, for any $w \in N_{r} (z)$, we have 
\begin{equation}
\tau_{kr} (z)\ \ \le\ \ \tau_{r} (w, z)\ \ \le\ \ \tau_{\frac{1}{k}r} (z)\ \ \ \ \forall k \ge 1. 
\end{equation}
\medskip

\noindent 
In this paper, we shall use the above definitions of the recurrence times to understand the concept of recurrence rate of typical orbits as $r \searrow 0$. We first define the \emph{recurrence rate of the centre} $z \in N_{r} (z)$ as follows. 
\begin{eqnarray*} 
\underline{R} (z) & := & - \liminf_{r \to 0} \frac{\log \tau_{r} (z)}{\log r} ; \\
\overline{R} (z) & := & - \limsup_{r \to 0} \frac{\log \tau_{r} (z)}{\log r}. 
\end{eqnarray*} 
If $\underline{R} (z) = \overline{R} (z)$, then the \emph{recurrence rate} is denoted and defined by, 
\begin{equation} 
R (z)\ \ :=\ \ - \lim_{r \to 0} \frac{\log \tau_{r} (z)}{\log r}. 
\end{equation} 
\medskip 

\noindent 
Now consider $\mu \in \mathcal{M}_{T} (\mathcal{J})$. We now define a local version of the \emph{fractal dimension} of $N_{r} (z)$ with respect to this chosen measure $\mu$ as follows. 
\begin{eqnarray*} 
\underline{d}_{\mu} (z) & := & \liminf_{r \to 0} \frac{\log \left( \mu \left( N_{r} (z) \right) \right)}{\log r} ; \\
\overline{d}_{\mu} (z) & := & \limsup_{r \to 0} \frac{\log \left( \mu \left( N_{r} (z) \right) \right)}{\log r}. 
\end{eqnarray*}
\medskip 

\noindent 
Now consider the real-valued H\"{o}lder continuous function $f$ defined on $\mathcal{J}$, given by $f = - s \log |T'|$. Then we know that $f$ could be suitably normalised in order that $\textrm{Pr} (- s \log |T'|) = 0$. In other words, there exists a unique $s \in \mathbb{R}$ such that $\textrm{Pr} (- s \log |T'|) = 0$. This unique value of $s$ is called the \emph{Hausdorff dimension} of the Julia set, $\mathcal{J}$. Furthermore, observe that the essential supremum of $\underline{d}_{\mu}$ is nothing but the Hausdorff dimension of $\mathcal{J}$, given by $s$. 
\bigskip 

\section{Main Results} 

We state the main results of this paper in this section. The proofs of the results are given in the following sections. 
\bigskip 

\noindent 
\begin{theorem}
\label{nt1} 
Let $T$ be a hyperbolic rational map restricted on its Julia set, $\mathcal{J}$. Let $N_{r} (z)$ be a neighbourhood of radius $r > 0$ about the point $z$ in $\mathcal{J}$ and let $\mu \in \mathcal{M}_{T} (\mathcal{J})$. Then for $\mu$-a.e. $z \in \mathcal{J}$, we have  
\begin{eqnarray*} 
1. & & \underline{R} (z)\ \ \le\ \ \underline{d}_{\mu} (z) ; \\ 
2. & & \overline{R} (z)\ \ \le\ \ \overline{d}_{\mu} (z).
\end{eqnarray*} 
\end{theorem} 
\bigskip

\noindent 
\begin{theorem}
\label{nt2} 
Let $T$ be a hyperbolic rational map restricted on its Julia set, $\mathcal{J}$. Let $Cov_{T} (f, g)$ decay at a super-polynomial rate $(f, g \in L^{2} (\mu))$. Let $N_{r} (z)$ be a neighbourhood of radius $r > 0$ about the point $z$ in $\mathcal{J}$ and let $\mu \in \mathcal{M}_{T} (\mathcal{J})$. Then for $\mu$-a.e. $z \in \mathcal{J}$, we have  
\begin{eqnarray*} 
1. & & \underline{R} (z)\ \ \ge\ \ \underline{d}_{\mu} (z) ; \\ 
2. & & \overline{R} (z)\ \ \ge\ \ \overline{d}_{\mu} (z).
\end{eqnarray*} 
\end{theorem} 
\bigskip

\noindent 
The following is an immediate corollary that follows from the statements of theorems \ref{nt1} and \ref{nt2}.
\bigskip 

\noindent 
\begin{corollary}
\label{nc1} 
Let $T$ be a hyperbolic rational map restricted on its Julia set, $\mathcal{J}$. Let $Cov_{T} (f, g)$ decay at a super-polynomial rate $(f, g \in L^{2} (\mu))$. Let $N_{r} (z)$ be a neighbourhood of radius $r > 0$ about the point $z$ in $\mathcal{J}$ and let $\mu \in \mathcal{M}_{T} (\mathcal{J})$. Then for $\mu$-a.e. $z \in \mathcal{J}$, we have  
\begin{eqnarray*} 
1. & & \underline{R} (z)\ \ =\ \ \underline{d}_{\mu} (z) ; \\ 
2. & & \overline{R} (z)\ \ =\ \ \overline{d}_{\mu} (z).
\end{eqnarray*} 
\end{corollary} 
\bigskip

\section{Proof of Theorem \ref{nt1}} 

We begin this section with two definitions; diametrically regular measures and weakly diametrically regular measures, as can be found in Federer, \cite{hf:69}. 
\medskip 

\noindent 
A measure $\mu$ is called \emph{diametrically regular} if there exists $k > 1$ and $c > 0$ such that 
\begin{equation} 
\mu \left( N_{kr} (z) \right)\ \ \le\ \ c \mu \left( N_{r} (z) \right),\ \ \ \ \forall z \in \mathcal{J},\ \ \ \ \text{and}\ \ \ \forall r > 0. 
\end{equation} 
\medskip 

\noindent 
A measure $\mu$ is called \emph{weakly diametrically regular} on a set $B \subset \mathcal{J}$ if there exists $k > 1$ such that for $\mu$-a.e. $z \in B$ and $\alpha > 0$, there exists $\delta > 0$ that satisfies, 
\begin{equation} 
\mu \left( N_{kr} (z) \right)\ \ \le\ \ \frac{1}{r^{\alpha}} \mu \left( N_{r} (z) \right)\ \ \ \text{whenever}\ \ \ r < \delta. 
\end{equation} 
\medskip 

\noindent 
Though the next lemma is obvious by the nomenclatures of the above defined two terms, it is imperative that we specify the constants that relate them. 
\bigskip 

\noindent 
\begin{lemma} 
\label{nl1}
Diametrically regular measures are weakly diametrically regular on $\mathcal{J}$. 
\end{lemma} 
\bigskip 

\noindent 
A $T$-invariant probability measure $\mu$ supported on $\mathcal{J}$ is weakly diametrically regular if for every fixed constant $k > 1$, there exists a $\delta \equiv \delta(z, \alpha) > 0$ so that for $\mu$-a.e. $z \in B \subset \mathcal{J}$ and every $\alpha > 0$, we have 
\begin{equation}
\mu \left( N_{kr} (z) \right)\ \ \le\ \ \frac{1}{r^{\alpha}} \mu \left( N_{r} (z) \right)\ \ \ \text{whenever}\ \ \ r < \delta. 
\end{equation}
\bigskip 

\noindent 
\begin{lemma}
\label{nl2} 
Any Borel probability measure on $\mathbb{C}$ is weakly diametrically regular on its support. 
\end{lemma} 
\bigskip 

\noindent 
\begin{proof} 
In order to prove the statement in lemma \ref{nl2}, we should show: For $\mu$-a.e. $z \in \mathbb{C}$, 
\begin{equation}
\mu \left( N_{\frac{1}{2^{n}}} (z) \right)\ \ \le\ \ n^{2} \mu \left( N_{\frac{1}{2^{n + 1}}} (z) \right), \end{equation} 
for sufficiently large $n \in \mathbb{Z}_{+}$. 
\medskip 

\noindent 
For $n \in \mathbb{Z}_{+}$ and $\delta > 0$, define 
\[ K_{n} (\delta)\ \ :=\ \ \left\{ z \in \mathcal{J} : \mu \left( N_{\frac{1}{2^{n + 1}}} (z) \right) < \delta \mu \left( N_{\frac{1}{2^{n}}} (z) \right) \right\}. \] 
Let $E \subset K_{n} (\delta)$ be a maximal $\frac{1}{2^{n + 2}}$-separated set. Then, 
\begin{eqnarray*} 
\mu \left(K_{n} (\delta)\right) & \le & \sum_{z \in E} \mu \left( N_{\frac{1}{2^{n + 1}}} (z) \right) \\ 
& \le & \sum_{z \in E} \delta \mu \left( N_{\frac{1}{2^{n}}} (z) \right). 
\end{eqnarray*} 
Observe that $E$ can be written as a finite union of $\frac{1}{2^{n}}$-separated sets, i.e., $E = \cup_{i = 1}^{M} E_{i}$ such that each $E_{i}$ is $\frac{1}{2^{n}}$-separated. Here, $M$ depends on $n$. Thus the sets, $\left\{ N_{\frac{1}{2^{n}}} (z_{i}) \right\}_{z_{i} \in E_{i}}$ are pairwise disjoint. Hence, 
\begin{eqnarray*} 
\mu \left(K_{n} (\delta)\right) & \le & \sum_{z \in E} \delta \mu \left( N_{\frac{1}{2^{n}}} (z) \right) \\ 
& \le & M \delta. 
\end{eqnarray*} 
Put $\delta = \frac{1}{n^{2}}$. Then we have obtained, 
\[ \mu \left(K_{n} \left(\frac{1}{n^{2}}\right)\right)\ \ \le\ \ \frac{M}{n^{2}}\ \ \ \ \forall n. \]
Thus, 
\begin{equation} 
\sum_{n \ge 1} \mu \left(K_{n} \left(\frac{1}{n^{2}}\right)\right)\ \ \le\ \ M \sum_{n \ge 1} \frac{1}{n^{2}}\ \ < \infty. 
\end{equation} 
\bigskip 

\noindent 
\begin{lemma}[Borel - Cantelli Lemma] 
For a sequence $\{B_{n}\}$ in the $\sigma$-algebra of $(X, \mathcal{B}, \mu)$ that satisfies $\sum_{n \ge 1} \mu (B_{n}) < \infty$, we have $\mu \left( \limsup_{n \to \infty} B_{n} \right) = 0$. 
\end{lemma} 
\bigskip 

\noindent 
An application of Borel - Cantelli lemma then says, 
\[ \mu \left( \limsup_{n \to \infty} K_{n} \left( \frac{1}{n^{2}} \right) \right)\ \ =\ \ 0. \]
In other words, the set of all points that satisfy 
\[ \mu \left( N_{\frac{1}{2^{n + 1}}} (z) \right)\ \ <\ \ \frac{1}{n^{2}} \mu \left( N_{\frac{1}{2^{n}}} (z) \right) \]
is of measure zero, for sufficiently large $n \in \mathbb{Z}_{+}$, whence our claim. 
\end{proof} 
\bigskip 

\noindent 
\begin{proof}(of Theorem \ref{nt1})
Consider the function $\delta (z, \cdot)$ in the definition of a weakly diametrically regular measure. Observe that for every fixed $z \in B \subset \mathcal{J},\ \delta (z, \cdot)$ is a measurable function. Fix $\alpha > 0$ and choose $\rho > 0$ such that 
\[ \mu (B) - \mu (G)\ \ \le\ \ \epsilon, \]
where $G = \left\{ z \in B \subset \mathcal{J} : \delta (z, \alpha) > \rho \right\}$. 
\medskip 

\noindent 
For any $r > 0,\ \lambda > 0$ and $z \in \mathcal{J}$, consider the set 
\[ A_{4r} (z)\ \ :=\ \ \left\{ w \in N_{4r} (z) : \tau_{4r} (w, z) \ge \frac{1}{\lambda} \frac{1}{\mu \left( N_{4r} (z) \right)} \right\}. \]
\medskip 

\noindent 
The following is the well-known Chebyschev's inequality. 
\bigskip 

\noindent 
\begin{theorem}[Chebyschev's inequality] 
Let $(X, \mathcal{B}, \mu)$ be a probability measure space. Let $f$ be a real-valued measurable function defined on $X$. Then for any $t > 0$, we have 
\[ \mu \left( \left\{ x \in X : f(x) \ge t \right\} \right)\ \ \le\ \ \frac{1}{t} \int_{X} f d \mu. \]
\end{theorem} 
\bigskip 

\noindent 
Using Chebyschev's inequality on the set $A_{4r} (z)$, we have 
\begin{eqnarray*} 
\mu \left( A_{4r} (z) \right) & \le & \lambda \mu \left( N_{4r} (z) \right) \int_{N_{4r} (z)} \tau_{4r} (w, z) d \mu (w) \\ 
& = & \lambda \mu \left( N_{4r} (z) \right) \mu \left( \left\{ w \in \mathcal{J} : \tau_{4r} (w, z) < \infty \right\} \right) \\ 
& \le & \lambda \mu \left( N_{4r} (z) \right). 
\end{eqnarray*} 
Since $N_{2r} (z) \subset N_{4r} (z)$, we have 
\[ \mu \left( A_{2r} (z) := \left\{ w \in N_{2r} (z) : \tau_{4r} (w, z) \ge \frac{1}{\lambda} \frac{1}{\mu \left( N_{4r} (z) \right)} \right\} \right)\ \ \le\ \ \lambda \mu \left(N_{4r} (z) \right). \]
\medskip 

\noindent 
Moreover, for $w \in N_{2r} (z)$, we have 
\[ \tau_{8r} (w) \mu \left( N_{2r} (w) \right)\ \ \le\ \ \tau_{4r} (w, z) \mu \left( N_{4r} (z) \right). \]
Hence, 
\[ \mu \left( \left\{ w \in N_{2r} (z) : \tau_{8r} (w) \mu \left( N_{2r} (w) \right) \ge \frac{1}{\lambda} \right\} \right)\ \ \le\ \ \lambda \mu \left( N_{4r} (z) \right). \]
\medskip 

\noindent 
Now we state a lemma that is useful to complete the proof of theorem \ref{nt1}. The proof of the lemma will be taken up after we complete the proof of the theorem. 
\bigskip 

\noindent 
\begin{lemma} 
\label{nl3} 
Let $\mu \in \mathcal{M}_{T} (\mathcal{J})$. Let $E \subset \mathcal{J}$ be a measurable set. Given $r > 0$, there exists a countable set $K \subset E$ such that 
\begin{enumerate} 
\item $N_{r} (z) \cap N_{r} (w) = \varphi$, for any two distinct points $z, w \in K$. 
\item $\mu \left( E \backslash \cup_{z \in K} N_{2r} (z) \right) = 0$. 
\end{enumerate} 
\end{lemma} 
\bigskip 

\noindent 
Define a quantity $D_{\alpha} (r)$ as, 
\begin{equation} 
D_{\alpha} (r)\ \ :=\ \ \mu \left( \left\{ w \in E : \tau_{8r} (w) \mu \left( N_{2r} (w) \right) \ge \left( \frac{1}{r} \right)^{2 \alpha} \right\} \right). 
\end{equation} 
Then observe that one can obtain an upper bound for $D_{\alpha} (r)$ as follows. 
\begin{eqnarray} 
\label{ne1}
D_{\alpha} (r) & = & \mu \left( \left\{ w \in E : \tau_{8r} (w) \mu \left( N_{2r} (w) \right) \ge \left( \frac{1}{r} \right)^{2 \alpha} \right\} \right) \nonumber \\ 
& \le & \sum_{z \in K} \mu \left( \left\{ w \in N_{2r} (z) : \tau_{8r} (w) \mu \left( N_{2r} (w) \right) \ge \left( \frac{1}{r} \right)^{2 \alpha} \right\} \right) \nonumber \\ 
& \le & r^{2 \alpha} \sum_{z \in K} \mu \left( N_{4r} (z) \right) \nonumber \\ 
& \le & \frac{r^{2 \alpha}}{r^{\alpha}} \sum_{z \in K} \mu \left( N_{r} (z) \right) \nonumber \\ 
& \le & r^{\alpha}. 
\end{eqnarray} 
\medskip 

\noindent 
If we choose $r$ to vanish at an exponential rate, i.e., $r = e^{-n}$, then by the inequality in (\ref{ne1}), we have 
\[ D_{\alpha} \left( \frac{1}{e^{n}} \right)\ \ \le\ \ \frac{1}{e^{n \alpha}}. \]
Further, 
\[ \sum_{n} D_{\alpha} \left( \frac{1}{e^{n}} \right)\ \ \le\ \ \sum_{n} \frac{1}{e^{n \alpha}}\ \ <\ \ \infty, \]
as $n$ grows larger. 
\medskip 

\noindent 
We complete the proof of theorem \ref{nt1} by invoking the Borel - Cantelli lemma again that asserts 
\[ \mu \left( \limsup_{n \to \infty} D_{\alpha} \left( \frac{1}{e^{n}} \right) \right)\ \ =\ \ 0. \]
In other words, for $\mu$-a.e. $z \in E$, we must have for sufficiently large $n$, 
\[ \tau_{8r} (z) \mu \left( N_{2r} (z) \right)\ \ \le\ \ \left( \frac{1}{r} \right)^{2 \alpha},\ \ \ \text{where}\ \ \ r = \frac{1}{e^{n}}. \]
Taking logarithms, we have for sufficiently large $n$, 
\begin{equation} 
\log \tau_{\frac{8}{e^{n}}} (z) + \log \mu \left( N_{\frac{2}{e^{n}}} (z) \right)\ \ \le\ \ 2 n \alpha. 
\end{equation} 
Thus, we have 
\begin{equation} 
\frac{1}{n} \log \tau_{\frac{8}{e^{n}}} (z)\ \ \le\ \ 2 \alpha - \frac{1}{n} \log \mu \left( N_{\frac{2}{e^{n}}} (z) \right), 
\end{equation} 
for sufficiently large $n$. 
\end{proof} 
\bigskip 

\noindent 
We now complete this section by writing the proof of lemma \ref{nl3}. This is only an elementary exercise in basic set theory. 
\bigskip 

\noindent 
\begin{proof}(of Lemma \ref{nl3})
Fix $z \in E$ and consider the family of subsets of $E$ around $z$, ordered by inclusion; 
\[ \mathcal{F}_{z}\ \ :=\ \ \left\{ N_{r} (z) \subset E\ \ \text{for various values of}\ \ r \right\}. \] 
Observe that $\mathcal{F}_{z}$ is a totally ordered set with a maximal element in $\mathcal{F}_{z}$. However, $E$ is only a partially ordered set; i.e., given $r > 0$, we can find a $w \in E$ such that $N_{r} (z) \cap N_{r} (w) = \phi$. Put all such $z$ and $w$ in $K$. Then by construction, $K$ is countable. We must still verify the second property that the lemma asserts. 
\medskip 

\noindent 
Consider the totally ordered family $\mathcal{F}_{z}$ of subsets of $E$ for every $z \in K$. By Zorn's lemma, there exists a maximal element for every totally ordered chain in $E$. Observe that with this maximal element, the second property is satisfied. 
\end{proof} 
\bigskip 

\section{Proof of Theorem \ref{nt2}}

\noindent 
We begin this section with the statement of a proposition. The motivation behind the proposition, given after its statement, clinches the proof of theorem \ref{nt2}.  
\bigskip 

\noindent 
\begin{proposition} 
\label{nl4} 
Let $T$ be a hyperbolic rational map restricted on its Julia set, $\mathcal{J}$. Let $Cov_{T} (f, g)\ (f, g \in L^{2} (\mu))$ decay at a super-polynomial rate. For $a > 0$, consider the set $\{ z \in \mathcal{J} : \underline{d}_{\mu} (z) \ge a \}$. Given $\delta,\ \epsilon > 0$, there exists $\rho > 0$ (depending on $z$) such that 
\begin{equation} 
T^{n} (z)\ \ \notin\ \ N_{r} (z)\ \ \ \text{for any}\ \ r \in (0, \rho)\ \ \ \text{and}\ \ n\ \in\ \mathbb{Z}_{+} \cap \left[ \frac{1}{r^{\delta}}, \frac{1}{(\mu (N_{r} (z)))^{1 - \epsilon}} \right]. 
\end{equation} 
\end{proposition} 
\bigskip 

\noindent 
We shall now briefly look at the motivation behind this proposition. A rigorous proof of the same is written in subsection 6.1. The definition of mixing implies $\mu (N \cap T^{-n} N) \to (\mu (N))^{2}$ as $n \to \infty$. Therefore, for large $n$, we have $\mu (N \cap T^{-n} N) \le 2 (\mu (N))^{2}$. Extending this line of argument, one may observe that 
\begin{equation} 
\mu \left( N \cap T^{-n} N \cap T^{-n - 1} N \cap \cdots T^{-n - l} N \right)\ \ \le\ \ 2 l (\mu (N))^{2}. 
\end{equation} 
Now suppose $l \le (\mu (N))^{\epsilon - 1}$, then 
\begin{equation} 
\label{ne2}
\mu \left( N \cap T^{-n} N \cap T^{-n - 1} N \cap \cdots T^{-n - l} N \right)\ \ \le\ \ 2 (\mu (N))^{\epsilon}. 
\end{equation} 
Making use of the fact that the covariance decays at a super-polynomial rate, one can then estimate the size of the connected neighbourhood $N \subset \mathcal{J}$. One may then use the Borel - Cantelli lemma to show that typical points exhibit this property. An appropriate condition on $\underline{R}$ should then help in completing the proof of the theorem. 
\bigskip 

\noindent 
\begin{lemma}
\label{nl4.1} 
Let $T$ be a hyperbolic rational map restricted on its Julia set, $\mathcal{J}$. Let $Cov_{T} (f, g)$ decay at a super-polynomial rate locally $(f, g \in L^{2} (\mu))$ . Then 
\begin{eqnarray*} 
1. & & \underline{R} (z)\ \ \ge\ \ \underline{d}_{\mu} (z) ; \\ 
2. & & \overline{R} (z)\ \ \ge\ \ \overline{d}_{\mu} (z),
\end{eqnarray*} 
for $\mu$-a.e. $z \in \{ z \in \mathcal{J} : \underline{R} (z) > 0 \}$. 
\end{lemma}
\bigskip 

\noindent 
\begin{proof} 
Fix $a > 0$ and consider the set $\{ z \in \mathcal{J} : \underline{R} (z) > a \}$. Then, by theorem \ref{nt1}, we have that 
\[ \left\{ z \in \mathcal{J} : \underline{R} (z) > a \right\}\ \ \subset\ \ \left\{ z \in \mathcal{J} : \underline{d}_{\mu} (z) > a \right\}. \] 
By the definition of $\underline{R}$, we know that $r^{a} \tau_{r} (z) \ge 1$ for sufficiently small $r$ where $z \in \{ z \in \mathcal{J} : \underline{R} (z) > a \}$. 
\medskip 

\noindent 
Put $\delta = a$ in proposition \ref{nl4}. Then given $\epsilon > 0$, we have for $\mu$-a.e. $z \in \{ z \in \mathcal{J} : \underline{R} (z) > a \}$, 
\[ \tau_{r} (z)\ \ \ge\ \ \frac{1}{\left( \mu( N_{r} (z)) \right)^{1 - \epsilon}},\ \ \ \text{provided}\ r\ \text{is sufficiently small.} \]
Thus, 
\begin{equation}
\underline{R} (z)\ \ \ge\ \ (1 - \epsilon) \underline{d}_{\mu} (z). 
\end{equation} 
The arbitrariness of $\epsilon$ completes the proof. The other inequality can be obtained in a similar fashion. 
\end{proof} 
\bigskip 

\noindent 
We now explore the validity of the set $\{ z \in \mathcal{J} : \underline{R} (z) > a \}$, i.e., does there exist a subset of $\mathcal{J}$ whose elements have a strictly positive recurrence rate. A definition and a lemma that is useful for this purpose has been studied by Ornstein and Weiss \cite{ow:93}. 
\medskip 

\noindent 
Let $\mathfrak{P}$ be a partition of $\mathcal{J}$ that has finitely many subsets $\{ N_{i} \}_{i = 1}^{M}$ such that 
\[ \bigcup_{i = 1}^{M} \overline{N_{i}}\ \ =\ \ \mathcal{J}\ ;\ \ \ \ \mu \left( N_{i} \cap N_{j} \right)\ \ =\ \ 0\ \ \ \text{whenever}\ \ i \ne j. \] 
By $\mathfrak{P} (z)$, we denote that element $N_{i}$ in the partition $\mathfrak{P}$ that contains $z$. Further, a dynamically refined partition is defined as,  
\[ \mathfrak{P}^{n}\ \ :=\ \ \mathfrak{P}\ \vee\ T^{-1} \mathfrak{P}\ \vee\ \cdots\ \vee\ T^{- (n - 1)} \mathfrak{P}. \]
\medskip 

\noindent 
Suppose for $\mu$-a.e. $z \in \mathcal{J}$, there exists a positive real number $\lambda$ (dependent on $z$) such that $N_{\frac{1}{e^{n \lambda}}} (z) \subset \mathfrak{P}^{n} (z)$ for all $n$ sufficiently large, then we say that the partition $\mathfrak{P}$ has a large interior. 
\bigskip 

\noindent 
\begin{theorem}\cite{ow:93} 
\label{owt}
Let $\mathfrak{P}$ be a partition with large interior. Then $\underline{R} (z) > 0$ for $\mu$-a.e. $z \in \mathcal{J}$. 
\end{theorem}
\bigskip 

\noindent 
Thus, by theorem \ref{owt} due to Ornstein and Weiss, the existence of the set $\{ z \in \mathcal{J} : \underline{R} (z) > 0 \}$ depends on the existence of a partition $\mathfrak{P}$ of $\mathcal{J}$ with a large interior. The next proposition sheds light on the existence of such partitions. 
\bigskip 

\noindent 
\begin{proposition} 
\label{np5} 
Let $T$ be a hyperbolic rational map restricted on its Julia set, $\mathcal{J}$. If for $\mu$-a.e. $z \in \mathcal{J}$ there exists positive constants $\alpha,\ \beta$ such that $T^{n}$ behaves like a Lipschitz function with Lipschitz constant $e^{n \alpha}$ on $N_{\frac{1}{n \beta}} (z)$ for all $n$ sufficiently large, then there exists a partition $\mathfrak{P}$ of $\mathcal{J}$ that has a large interior. 
\end{proposition} 
\bigskip 

\noindent 
We now make use of proposition \ref{np5} to prove theorem \ref{nt2}. The proof of this proposition shall be given in subsection 6.2. We first prove that hyperbolic rational maps restricted on their Julia sets satisfy the hypothesis of proposition \ref{np5} and will then complete the proof of theorem \ref{nt2}.  
\bigskip 

\noindent 
\begin{proof}(of Theorem \ref{nt2}) 
As earlier, let $\mathfrak{P}$ be a partition of $\mathcal{J}$ into finitely many subsets $\{ N_{i} \}$ such that 
\[ \bigcup_{i = 1}^{M} \overline{N_{i}}\ \ =\ \ \mathcal{J}\ ;\ \ \ \ \mu \left( N_{i} \cap N_{j} \right)\ \ =\ \ 0\ \ \ \text{whenever}\ \ i \ne j. \] 
Observe that there exists a positive constant $\kappa (N_{i})$ such that 
\[ \left| T z_{1} - T z_{2} \right|\ \ \le\ \ \kappa (N_{i}) \left| z_{1} - z_{2} \right|\ \ \ \forall z_{1}, z_{2} \in N_{i}. \]
Define 
\[ \log K\ \ :=\ \ \int \log^{+} \kappa \left( \mathfrak{P} (z) \right) d \mu (z)\ \ =\ \ \sum_{N_{i} \in \mathfrak{P}} \log^{+} \kappa (N_{i}) \mu (N_{i}). \]
Now choose $\alpha > \log K$. Observe that by the definition of ergodicity, as stated in section 2 the following statement is true. For $\mu$-a.e. $z \in \mathcal{J}$, there exists a $m \in \mathbb{Z}_{+}$ (dependent on $z$) that satisfies 
\begin{equation} 
\label{ne3}
\kappa \left( \mathfrak{P} (z) \right) \times \kappa \left( \mathfrak{P} (Tz) \right) \times \cdots \times \kappa \left( \mathfrak{P} (T^{n - 1}z) \right)\ \ \le\ \ e^{n \alpha}\ \ \ \forall n \ge m. 
\end{equation} 
Replace the upper bound in inequality \ref{ne3} by $c(z) e^{n \alpha}$ for some constant $0 \le c(z) \le 1$ in order that the above inequality is valid for every positive integer $n \in \mathbb{Z}_{+}$. 
\medskip 

\noindent 
Now choose $\beta > 0$ such that 
\[ N_{\frac{1}{c(z) e^{\beta}}} (Tz)\ \ \subset\ \ \mathfrak{P} (Tz). \] 
It is then a simple exercise to prove (by induction on $n$) that 
\begin{equation} 
N_{\frac{1}{c(z) e^{n \beta}}} (T^{n} z)\ \ \subset\ \ \mathfrak{P} (T^{n} z).
\end{equation} 
\medskip 

\noindent 
Thus, we can apply proposition \ref{np5} to say that there exists a partition $\mathfrak{P}$\ of $\mathcal{J}$ that has a large interior. A result due to Ornstein and Weiss as stated in theorem \ref{owt} and lemma \ref{nl4.1} then completes the proof of theorem \ref{nt2}. 
\end{proof} 
\bigskip 

\subsection{Proof of Proposition \ref{nl4}}

\noindent 
Consider the real-valued H\"{o}lder continuous function $f = - s \log |T'|$ defined on $\mathcal{J}$. Then we know that there exists a unique value of $s \in \mathbb{R}$ called the Hausdorff dimension of $\mathcal{J}$ so that $\textrm{Pr} (-s \log |T'|) = 0$. Let $a > 0$ be as given in the statement of proposition \ref{nl4}. Fix $b > 0$ and let $c = \frac{a}{3} \epsilon$ for some $\epsilon > 0$. Denote by $\mathcal{J}_{a}$ the set, 
\[ \mathcal{J}_{a}\ \ :=\ \ \left\{ z \in \mathcal{J} : \underline{d}_{\mu} (z) \ge a \right\}. \] 
For some $\rho > 0$, consider the following sets. 
\begin{eqnarray*} 
G_{1} & := & \left\{ z \in \mathcal{J}_{a} : \mu \left( N_{r} (z) \right) \le r^{a}\ \ \forall r \le \rho \right\} \\ 
G_{2} & := & \left\{ z \in \mathcal{J} : \mu \left( N_{r} (z) \right) \ge r^{b + s}\ \ \forall r \le \rho \right\} \\ 
G_{3} & := & \left\{ z \in \mathcal{J} : \mu \left( N_{\frac{r}{2}} (z) \right) \ge r^{c} \mu \left( N_{4r} (z) \right)\ \ \forall r \le \rho \right\} 
\end{eqnarray*}
Observe that by the definition of lower point wise dimension, we have that $\mu (G_{1}) \to \mu (\mathcal{J}_{a})$ as $\rho \to 0$. Further, $\mu (G_{i}) \to 1$ as $\rho \to 0$ for $i = 2, 3$, the former since $\overline{d}_{\mu} (z) \le s$ a.e. while the latter because the measure $\mu$ is diametrically regular. Thus, defining $G := G_{1} \cap G_{2} \cap G_{3}$, we claim $\mu (G) \to \mu (\mathcal{J}_{a})$. 
\medskip 

\noindent 
For $r \le \rho$, define 
\begin{equation}
B_{\epsilon} (r)\ \ :=\ \ \left\{ z \in \mathcal{J} : \exists n \in \mathbb{Z}_{+} \cap \left[ \frac{1}{r^{\delta}}, \frac{1}{(\mu(N_{3r} (z)))^{1 - \epsilon}} \right]\ \ \text{satisfying}\ \ n \ge \tau_{r} (z) \right\} 
\end{equation}
For $z \in \mathcal{J}$, observe  
\begin{eqnarray*} 
N_{r} (z) \cap B_{\epsilon} (r) & \subset & \left\{ w \in N_{r} (z) : \exists n \in \mathbb{Z}_{+} \cap \left[ \frac{1}{r^{\delta}}, \frac{1}{(\mu(N_{2r} (z)))^{1 - \epsilon}} \right]\ \ \ \text{satisfying}\ \ n \ge \tau_{2r} (w, z) \right\} \\ 
\vspace{+10pt}
& \subset & \bigcup_{n \in \mathbb{Z}_{+} \cap \left[ \frac{1}{r^{\delta}}, \frac{1}{(\mu(N_{2r} (z)))^{1 - \epsilon}} \right]} N_{r} (z) \cap T^{-n} N_{2r} (z). 
\end{eqnarray*}
Let $\varphi_{r} : [0, \infty) \longrightarrow \mathbb{R}$ be a Lipschitz function with Lipschitz constant $r^{-1}$ such that it is sandwiched between the characteristic functions of the intervals $[0, r]$ and $[0, 2r]$. Fix $z \in \mathcal{J}$ and define $\psi_{z, r} : \mathcal{J} \longrightarrow \mathbb{R}$ by $\psi_{z, r} (w) := \varphi_{r} ( |z - w| )$. It is a simple observation that $\psi_{z, r}$ is a Lipschitz function with the same Lipschitz constant as $\varphi_{r}$. 
\medskip 

\noindent 
\begin{eqnarray*} 
\mu \left( N_{r} (z) \cap T^{-n} N_{2r} (z) \right) & \le & \int \psi_{z, 2r} \psi_{z, 2r} \circ T^{n} d \mu \\ 
& \le & \| \psi_{z, 2r} \|^{2} \theta_{n}\ +\ \left( \int \psi_{z, 2r} d \mu \right)^{2} \\ 
& \le & \frac{\theta_{n}}{r^{2}}\ +\ \left( \mu \left( N_{4r} (z) \right) \right)^{2}. 
\end{eqnarray*} 
Here, to obtain the second inequality, we have used the hypothesis of lemma \ref{nl4} that the covariance decays at a super-polynomial rate and to obtain the third inequality, we have used the fact that $\psi_{z, r}$ is a Lipschitz function. 
\medskip 

\noindent 
Choose $p > 1$ sufficiently large such that $\delta (p - 1) \ge s + 2b + 2$. Further, choose $\rho > 0$ so small in order that 
\[ \frac{1}{\rho^{\delta}}\ \ \le\ \ n\ \ \ \ \Longrightarrow\ \ \ \ \theta_{n}\ \ \le\ \ \frac{p - 1}{n^{p}}. \]
Since 
\[ \sum_{n \ge q} \frac{1}{n^{p}}\ \ \le\ \ \frac{1}{p - 1} \frac{1}{q^{p - 1}}, \]
we obtain 
\begin{eqnarray} 
\mu \left( N_{r} (z) \cap B_{\epsilon} (r) \right) & \le & r^{\delta (p - 1) - 2} + \frac{(\mu(N_{4r} (z)))^{2}}{(\mu(N_{2r} (z)))^{1 - \epsilon}} \nonumber \\ 
& \le & \mu(N_{\frac{r}{2}} (z)) \left(r^{b} + r^{a \epsilon - 2c}\right). 
\end{eqnarray} 
\medskip 

\noindent 
Let $B \subset G$ be a maximal $r$-separated set, i.e., $N_{r} (\zeta) \cap B = \{ \zeta \}$ for every $\zeta \in B$ and for every $\omega \in G \setminus B$ there exists a $\zeta \in B$ such that $\zeta \in N_{r} (\omega)$. In other words, the family of neighbourhoods $\{ N_{r} (\zeta) \}_{\zeta \in B}$ is an open cover for $G$. Thus, 
\begin{eqnarray} 
\mu \left( G \cap B_{\epsilon} (r) \right) & \le & \sum_{\zeta \in B} \mu \left( N_{r} (\zeta) \cap B_{\epsilon} (r) \right) \nonumber \\ 
& \le & \sum_{\zeta \in B} \mu(N_{\frac{r}{2}} (\zeta)) \left(r^{b} + r^{a \epsilon - 2c}\right) \nonumber \\ 
& \le & r^{b} + r^{a \epsilon - 2c}.
\end{eqnarray} 
\medskip 

\noindent 
Choosing $r$ to vanish at an exponential rate, i.e., $r = e^{-n}$, we obtain 
\begin{equation} 
\sum_{n \ge 0} \mu \left( B_{\epsilon} \left( \frac{1}{e^{n}} \right) \right)\ \ <\ \ \infty. 
\end{equation} 
\medskip 

\noindent 
Applying the Borel-Cantelli lemma, we then obtain that for $\mu$-a.e. $w \in G$ there exists a $m(w) \in \mathbb{Z}_{+}$ such that for every $m > m(w)$ there does not exist any 
\begin{equation} 
n \in \mathbb{Z}_{+} \cap \left[ \frac{1}{e^{m \delta}}, \frac{1}{\left(\mu (N_{\frac{3}{e^{m}}} (w))\right)^{1 - \epsilon}} \right]\ \ \ \text{such that}\ \ \ T^{n} w \in N_{\frac{1}{e^{m}}} (w). 
\end{equation}
\medskip 

\noindent 
The weak diametric regularity of $\mu$ then completes the proof of proposition \ref{nl4}. 
\bigskip 

\subsection{Proof of Proposition \ref{np5}} 

\noindent 
Before we embark on proving proposition \ref{np5}, we state and prove the following lemma. 
\bigskip 

\noindent 
\begin{lemma} 
\label{nl6}
For every $z \in \mathcal{J}$, and for any $r > 0$, there exists a $\rho \in (r, 2r)$ such that 
\[ \mu \left( \left\{ w \in \mathcal{J} : \rho - \frac{r}{4^{n + 1}} < | z - w | < \rho + \frac{r}{4^{n + 1}} \right\} \right)\ \ \le\ \ \frac{1}{2^{n}} \mu \left( N_{2r} (z) \right). \]
\end{lemma}
\bigskip 

\noindent 
\begin{proof} 
Define a measure $m$ on the interval $(0, 2)$ by $m( [0, t) ) := \mu( N_{rt} (z) )$. Let $I_{0} = (1, 2)$. Divide $I_{0}$ into four pieces of equal length and define $I_{1}$ to be that interval which satisfies, $m(I_{1}) \le m(I_{0})/2$. Proceeding on, define $I_{n + 1}$ to that interval which satisfies, $m(I_{n + 1}) \le m(I_{n})/2$. Thus, we have constructed a decreasing sequence of nested intervals $\{I_{n}\}_{n \ge 0}$, in whose intersection must lie a single point, say $R$. Thus, 
\[ m \left( \left( R - \frac{1}{4^{n + 1}}, R+ \frac{1}{4^{n + 1}} \right) \right)\ \ \le\ \ m\left(I_{n}\right)\ \ \le\ \ \frac{1}{2^{n}} m\left(I_{0}\right). \] 
\end{proof} 
\bigskip 

\noindent 
We now complete the proof of proposition \ref{np5} and thus the proof of theorem \ref{nt2}. 
\bigskip 

\noindent 
\begin{proof}(of proposition \ref{np5}) 
Fix $r > 0$. Let $\mathfrak{P} = \{ N_{r} (z) \}$ be some partition of $\mathcal{J}$. Choose a maximal $r$-separated set $E \in \mathfrak{P}$. For any $z \in E$, take $\rho_{z} \in (r, 2r)$ such that lemma \ref{nl6} holds. Let $E = \{z_{1}, z_{2}, \cdots \}$ be an enumeration of the (at most) countable set $E$. Write $N_{i} = N_{\rho_{z_{i}}} (z_{i})$ and define 
\[ Q_{i}\ \ :=\ \ N_{i} \setminus \left( Q_{1} \cup Q_{2} \cup \cdots \cup Q_{i - 1} \right)\ \ \ \ \text{starting with}\ \ \ \ Q_{1}\ \ =\ \ N_{1}. \] 
\medskip 

\noindent 
Observe that by the maximality of the collection of sets, $\mathfrak{Q} := \{ Q_{i} \}_{i \ge 1}$ is a partition of $\mathcal{J}$. Further, $\partial \mathfrak{Q} \subset \bigcup_{i} \partial N_{i}$. Thus, 
\begin{eqnarray} 
\mu \left( \left\{ z \in \mathcal{J} : d(z, \partial \mathfrak{Q}) < \frac{r}{4^{n + 1}} \right\} \right) & \le & \mu \left( \bigcup_{i} \left\{ z \in \mathcal{J} : \rho_{z_{i}} - \frac{1}{4^{n + 1}} < |z - z_{i}| < \rho_{z_{i}} + \frac{1}{4^{n + 1}} \right\} \right) \nonumber \\ 
& \le & \frac{1}{2^{n}} \sum_{i} \mu \left( N_{2r} (z_{i}) \right). 
\end{eqnarray} 
\medskip 

\noindent 
Recall that $z_{i}$ were so chosen that they were $r$-separated. Hence, there exists at most $c = c (\dim (\mathcal{J}) = s)$ balls of radius $2r$ that can intersect $\mathcal{J}$. Thus for some constants $a,\ c > 0$ and for all $\epsilon > 0$, we have 
\[ \mu \left( \left\{ z \in \mathcal{J} : d(z, \partial \mathfrak{Q}) < \epsilon \right\} \right)\ \ <\ \ c \epsilon^{a}. \] 
\medskip 

\noindent 
Thus for any $b > 0$ we have by the invariance of $\mu$, 
\[ \sum_{n} \mu \left( \left\{ z \in \mathcal{J} : d ( T^{n} z, \mathfrak{Q}) < \frac{1}{e^{bn}} \right\} \right)\ \ \le\ \ \sum_{n} \frac{c}{e^{abn}}\ \ <\ \ \infty. \]
\medskip 

\noindent 
Hence, by the Borel - Cantelli Lemma 
\begin{equation} 
\exists n(z) < \infty\ \ \ \ \text{such that}\ \ \ \ d(T^{n} z, \partial \mathfrak{Q}) \ge \frac{1}{e^{bn}},\ \ \ \ \text{for}\ \ \ \ \mu \text{-a.e.}\ z \in \mathcal{J}. 
\end{equation}
This implies, $N_{\frac{1}{e^{bn}}} (T^{n} z) \subset \mathfrak{Q} (T^{n} z)$, for any $n \ge n(z)$. We can always choose $c(z) \in (0, 1)$ sufficiently small in order that 
\begin{equation}
N_{\frac{c(z)}{e^{bn}}} (T^{n} z)\ \ \subset\ \ \mathfrak{Q} (T^{n} z)\ \ \ \ \forall n \in \mathbb{Z}_{+}. 
\end{equation}
\medskip 

\noindent 
Choose $z \in \mathcal{J}$ such that the hypothesis in proposition \ref{np5} holds, i.e., let $z \in \mathcal{J}$ be a point that assures the existence of positive constants $\alpha,\ \beta$ such that $T^{n}$ behaves like a Lipschitz function with Lipschitz constant $e^{n \alpha}$ on $N_{\frac{1}{n \beta}} (z)$ for $n$ sufficiently large. Observe that without of generality, we can change $c(z)$ into a smaller constant (if necessary) so that $T^{n}$ remains $e^{n \alpha}$-Lipschitz on $N_{\frac{c(z)}{e^{n \beta}}} (z)$ for all integers $n$ and that $\alpha > b + \beta$.
\medskip 

\noindent 
The proof is then complete if we prove 
\begin{equation} 
N_{\frac{c(z)}{e^{n \alpha}}} (z)\ \ \subset\ \ \mathfrak{Q}^{k} (z)\ \ \ \ \text{for any}\ \ \ \ k \le n. 
\end{equation}
We prove this by induction. Observe that the statement is trivially true for $k = 1$. Assuming the statement to be true for $k \le n - 1$, observe that 
\begin{equation} 
T^{k} \left( N_{\frac{c(z)^{2}}{e^{n \beta}}} (z) \right)\ \ \subset\ \ N_{\frac{c(z) e^{(k} \alpha}{n \beta}} (T^{k} z)\ \ \subset\ \ N_{\frac{1}{e^{b n}}} (T^{k} z)\ \ \subset \mathfrak{Q} (T^{k} z).
\end{equation} 
Thus, $N_{\frac{c(z)^{2}}{e^{n \beta}}} (z) \subset \mathfrak{Q}^{k + 1} (z)$. 
\end{proof} 
 
\bigskip 
\bigskip 
\bigskip 

\noindent 
\textbf{Shrihari Sridharan} \\
\bigskip

\noindent 
\textsc{Chennai Mathematical Institute (CMI)}, \\ 
Plot \# H1, SIPCOT IT Park, Kelambakkam, Siruseri, Chennai, India.\ \ \ \ PIN 603 103. \\ 
email: \texttt{shrihari@cmi.ac.in}


\bigskip \bigskip 

\begin{thebibliography}{88} 

\bibitem{bs:01} \textsc{Barreira, L.} and \textsc{Saussol, B.}, "Hausdorff dimension of measures via Poincar\'{e} recurrence", \emph{Communications in Mathematical Physics}, {\bf 219}, (2001) pp 443 - 463.

\bibitem{mb:93} \textsc{Boshernitzan, M. D.}, "Quantitative recurrence results", \emph{Inventiones Mathematicae}, {\bf 113}, (1993) pp 617 - 631. 

\bibitem{du:91} \textsc{Denker, M.} and \textsc{Urbanski, M.}, "Ergodic theory of equilibrium states for rational maps", \emph{Nonlinearity}, {\bf 4}, (1991), pp 103 - 134.

\bibitem{hf:69} \textsc{Federer, H.}, "Geometric measure theory", BerlinÐHeidelbergÐNew York, Springer, 1969.

\bibitem{nh:99} \textsc{Haydn, N.}, "Convergence of the transfer operator for rational maps", \emph{Ergodic Theory and Dynamical Systems}, {\bf 19}, (1999), pp 657 - 669.

\bibitem{ml:86} \textsc{Lyubich, M. Yu.}, "The dynamics of rational transforms: the topological picture", \emph{Russian Math. Surveys}, (1986), pp 43 - 117.

\bibitem{ow:93} \textsc{Ornstein, D. S.} and \textsc{Weiss, B.}, "Entropy and data compression schemes", \emph{IEEE Transactions on Information Theory}, {\bf 39}, (1993), pp 78 - 83. 

\bibitem{bs:06} \textsc{Saussol, B.}, "Recurrence rates in rapidly mixing dynamical systems", \emph{Discrete Contin. Dyn. Syst.}, {\bf 15}, (2006), pp 259 - 267. 

\bibitem{ss:12} \textsc{Sridharan, S.}, "Sinai - Ruelle - Bowen measure leaks", \emph{Communications in Mathematical Analysis}, {\bf 13}, (2012) pp 1 - 22.

\bibitem{pw:82} \textsc{Walters, P.}, "An introduction to ergodic theory", \emph{Graduate Texts in Mathematics}, Springer-Verlag, New York, 79, (1982).

\end{thebibliography}
\end{document}